
\input amstex.tex
\documentstyle{amsppt}

                                                 \newif\iffigs\figsfalse


\figstrue

\def\pretitle{\noindent\hskip23.5pc\vbox to0pt{%
                    \vskip-31pt\advance\hsize-23.5pc\parindent0pt
     DUK-M-94-04    \hfil\break
     June, 1994     \par\vss}\vskip3pc}

\magnification\magstep1
\pageheight{38pc}
\pagewidth{29pc}
\voffset3pc
\hoffset1pc

\iffigs
  \input epsf
\else
  \message{No figures will be included. See TeX file for more
information.}
\fi

\define\eqref#1{(#1)}

\define\GW{2.1}
\define\threepoint{2.2}
\define\sing{4.1}
\define\intheform{4.2}
\define\firstterms{4.3}
\define\balance{4.4}

\define\catp{1}
\define\small{2}
\define\AM{3}
\define\StringKThree{4}
\define\suffice{5}
\define\BatyrevDuke{6}
\define\Bogomolov{7}
\define\BurnsRapoport{8}
\define\Calabi{9}
\define\lagrangian{10}
\define\CdO{11}
\define\twoparam{12}
\define\twoparamtwo{13}
\define\CLS{14}
\define\CKM{15}
\define\Dixon{16}
\define\DH{17}
\define\Friedman{18}
\define\GlimmJaffe{19}
\define\GrassiMorrison{20}
\define\GP{21}
\define\Gromov{22}
\define\GS{23}
\define\Hirz{24}
\define\Katz{25}
\define\appendix{26}
\define\Kawamata{27}
\define\Kirwan{28}
\define\Kontsevich{29}
\define\KontsevichManin{30}
\define\LVW{31}
\define\LiuMcDuff{32}
\define\Looijenga{33}
\define\McDuff{34}
\define\remarksKthree{35}
\define\compact{36}
\define\Namikawa{37}
\define\Ness{38}
\define\Ruan{39}
\define\RuanTian{40}
\define\schoen{41}
\define\Thaddeus{42}
\define\Tian{43}
\define\Todorov{44}
\define\Vafa{45}
\define\Viehweg{46}
\define\Wilson{47}
\define\tsm{48}
\define\topgrav{49}
\define\Witten{50}
\define\phases{51}
\define\Yau{52}
\define\YauBook{53}

\let\em\it
\let\cal\Cal

\define\C{\Bbb C}
\define\cC{{\cal C}}
\redefine\D{{\cal D}}
\define\F{{\cal F}}

\define\K{{\cal K}}
\define\M{{\cal M}}
\define\IM{\Bbb M}
\redefine\O{{\cal O}}
\redefine\P{\Bbb P}
\define\Q{\Bbb Q}
\define\R{\Bbb R}
\define\T{{\cal T}}
\define\X{{\cal X}}
\define\Z{\Bbb Z}

\define\Aut{\operatorname{Aut}}
\define\Bir{\operatorname{Bir}}
\define\DR{\operatorname{DR}}
\redefine\Im{\operatorname{Im}}
\define\PSL{\operatorname{PSL}}
\define\Spec{\operatorname{Spec}}
\define\SU{\operatorname{SU}}
\define\U{\operatorname{U}}

\define\Mov{\mathop{{\cal M}ov}\nolimits}

\define\Ms{\M_\sigma}
\define\MN{\M_{N{=}2}}
\define\Mc{\M_{\text{complex}}}
\define\Mbir{\M_{N{=}2}^{\text{birat}}}

\define\bra{\langle}
\define\ket{\rangle}
\define\semidirect{\rtimes}

\define\catquot{\mathrel{/\!/}}

\topmatter
\title Beyond the K\"ahler cone \endtitle
\author David R. Morrison
\endauthor
\address Duke University and the Institute for Advanced Study \endaddress
\email drm\@math.duke.edu \endemail
\dedicatory To Friedrich Hirzebruch \enddedicatory
\thanks Research partially supported by  National Science Foundation
Grant DMS-9103827, and by
an American Mathematical Society Centennial Fellowship. \endthanks
\subjclass Primary 14J30; Secondary 14E07, 14J15, 32J81, 58D27,
81T40 \endsubjclass
\abstract
The moduli space of nonlinear $\sigma$-models on a Calabi--Yau manifold
contains a complexification of the K\"ahler cone of the manifold.  We
describe a physically natural analytic continuation process which links
the complexified K\"ahler cones of birationally equivalent Calabi--Yau
manifolds.  The enlarged moduli space includes a complexification of
Kawamata's ``movable cone''.  We formulate a natural conjecture about the
action of the birational automorphism group on this cone.
\endabstract
\endtopmatter

\document

Many mathematicians were taken by surprise during 1984--85 when we found
physicists knocking on our doors, asking whether we knew anything about
Riemannian 6-manifolds with a metric whose holonomy lies in $\SU(3)$.
Fortunately, Yau had solved the Calabi conjecture nearly 10 years earlier,
so we were able to provide some answers:  any smooth complex projective
threefold with trivial canonical bundle admits a metric of this type.
Also fortunately, these mani\-folds---now called Calabi--Yau
threefolds---had been studied in some detail
by algebraic geometers, in part due to the distinguished r\^ole they
play in the classification theory of algebraic varieties.

During the following year, the questions became more and more specific,
focusing primarily on the physicists' desire for examples $X$ whose Euler
number $e(X)$ satisfies $e(X)=\pm6$.\footnote{In retrospect, perhaps
someone should have wondered why the
physicists didn't know whether they wanted $e(X)=6$
or $e(X)=-6$.  This ambiguity
was one of the early hints of the phenomenon of ``mirror symmetry'',
which I will discuss later.}
A lot of work was done on this question at the Max-Planck-Institut
in Bonn during 1985--86, and much of it is described by Hirzebruch in
the notes of his very nice lectures \cite{\Hirz}.
(Other papers from the same period which explain the r\^ole of
Calabi--Yau manifolds in physics can be found in \cite{\YauBook};
 cf.\ also \cite{\Friedman}.)
One of my purposes  in this paper is to update the story, and explain what
has been happening recently in the study of Calabi--Yau manifolds from
the point of view of physics.

The focus of the investigations by string theorists has shifted from
a search for specific Calabi--Yau threefolds to a study of the
general properties of the physical theories built from such manifolds.
The subject was re-invigorated a few
years ago by the discovery of the surprising phenomenon known as mirror
symmetry.  I will review this, with a focus on recent applications of
 mirror symmetry  to the study of the moduli spaces of the physical
theories.  It has recently been discovered \cite{\phases, \catp} that
there is a kind of analytic continuation which links the physical theories
associated to different birational models of a single Calabi--Yau manifold.
The primary purpose of this paper is to explain this analytic continuation
in a very general setting, for arbitrary Calabi--Yau threefolds.
This work is an outgrowth of my collaboration
 with Paul Aspinwall and Brian Greene \cite{\catp, \small},
whose contribution I would like to acknowledge at the outset.

\head 1.~Moduli spaces of $\sigma$-models \endhead

A {\em Calabi--Yau manifold}\/ is a compact connected orientable
 manifold $X$ of dimension $2n$
which admits Riemannian metrics whose (global) holonomy lies in $\SU(n)$.
Physicists have constructed some two-dimensional
quantum field theories associated to these
manifolds which are known as
 {\em nonlinear $\sigma$-models}.  Quantum field theories
remain problematic for mathematicians, since they have not been shown
to make sense as rigorous mathematical theories except in very limited
cases (cf.~\cite{\GlimmJaffe}).  However, certain aspects of the quantum
field theories built from Calabi--Yau manifolds can be studied purely
mathematically, without making reference to the underlying physical
theory.  This is the strategy we shall adopt here.

The two-dimensional quantum field theories in question
are constructed out of the space of maps from (variable) surfaces
$\Sigma$ to a (fixed) Calabi--Yau manifold $X$.  Each such
theory is based on
 a ``Lagrangian''
functional on the space of such maps, and the standard Lagrangian
used (cf.~\cite{\lagrangian}) depends on the choice of a pair $(g_{ij},B)$,
where $g_{ij}$ is a Riemannian metric on $X$ with holonomy contained
in\footnote{More precisely, the metrics which are needed in
the physical theory are perturbations of these metrics with restricted
holonomy.  However the perturbed metrics, like the original metrics
with restricted holonomy, are expected to be uniquely determined by
the cohomology class of the associated K\"ahler form,
once a complex structure has been specified on $X$.}
$\SU(n)$,
and $B$ is the de Rham cohomology class of a real closed $2$-form on $X$.
We consider two such pairs to be equivalent, written
$(g_{ij},B)\sim(g_{ij}',B')$, if there is a diffeomorphism $\varphi:X\to X$
such that $\varphi^*(g_{ij})=g_{ij}'$ and
$\varphi^*(B)-B'\in H^2_{\DR}(X,\Z)$, where $H^2_{\DR}(X,\Z)$ denotes
the image of the integral cohomology in the de Rham cohomology.
This definition of equivalence arises as follows.  First, the appearance
of $B$ in the Lagrangian (when applied to the map $f:\Sigma\to X$) takes the
form $\int_\Sigma f^*(B)$, so only the de Rham class of $B$
matters.\footnote{Note that this term---and indeed
the entire Lagrangian---is
 also invariant under a simultaneous change of the sign of $B$ and
the orientation of $\Sigma$.  (I am grateful to Paul Aspinwall and
Jacques Distler for discussions on this point.)  We can safely ignore
this here, though, since we have not specified an orientation of $\Sigma$.}
Second, the appearance of this Lagrangian in physically measurable
quantities always involves an exponentiation in which this term
becomes $\exp(2\pi i\int_\Sigma f^*(B))$.  Thus, shifting $B$ by an
integral class will not affect the physical theory.\footnote{A bit
more generally, one should also include a contribution to
$\exp(2\pi i\int_\Sigma f^*(B))$ coming from torsion in $H_2(X,\Z)$,
as in \cite{\Vafa, \suffice};
we will suppress that contribution in this
paper.}

As in \cite{\compact}, we regard the set of equivalence classes of pairs
$$\Ms:=\{(g_{ij},B)\}/\sim$$
as a first approximation to a moduli space for these theories, which
we call the {\em one-loop semiclassical nonlinear $\sigma$-model moduli
space}, or just the {\em nonlinear $\sigma$-model moduli space}\/ for short.
This space may differ from the actual moduli space in several ways.
First, the physical theory may fail to converge for some values of
$(g_{ij},B)$.  This statement has no mathematical content in the absence
of an adequate mathematical definition of quantum field theories; however,
certain parts of the physical theory which can be formulated in purely
mathematical terms (such as the three-point functions described below)
should be expected to converge, and their failure to converge at a certain
place
is evidence that the physical theory is badly behaved there.

Second, the family of quantum field theories may admit an analytic
continuation (regarding them purely
as two-dimensional quantum field theories) which no longer
has a nonlinear $\sigma$-model interpretation.
In fact, our two-dimensional quantum field theories are of
a type called ``superconformal'', and we are  interested
in deformations which preserve only this ``superconformal'' property
and not the more restrictive
``$\sigma$-model'' property.  Varying
$g_{ij}$ and the class of $B$ gives a locally complete family of such
deformations, but globally there may be  deformations
which do not preserve the $\sigma$-model structure.

Third, two superconformal field theories may be isomorphic by an isomorphism
which does not preserve the $\sigma$-model structure.  So we may have to
enlarge the set of identifications among pairs which we are making.

In sum, we may need to shrink our moduli space a bit to ensure convergence,
we may need to enlarge it (in other directions) to get a complete family,
and we may need to mod out by further discrete identifications.  In spite
of these limitations,
we can still obtain
 significant information about the structure of the moduli space
by studying the more primitive ``one-loop semiclassical
nonlinear $\sigma$-model'' version we
have formulated above.

We immediately need a slight refinement of the nonlinear $\sigma$-model moduli
space, which we call the {\em $N{=}2$ moduli space}.
This is defined to be
$$\MN:=\{(g_{ij},B,t)\}/\sim$$
where $t$ denotes a complex structure on $X$ with respect to which
$g_{ij}$ is K\"ahler.
(That such complex structures exist is a consequence of the holonomy
being contained in $\U(n)$.)
By the Bogomolov--Tian--Todorov theorem \cite{\Bogomolov, \Tian, \Todorov},
 deformations
of complex structure on $X$ are unobstructed and  the Kodaira--Spencer
map of each versal deformation  is an isomorphism.
By way of notation, we let $\X_t$ denote
the complex manifold ``$X$ equipped with the structure $t$'', which
we sometimes treat as if it were a fiber of a universal family
$\X\to\Mc$ over the moduli space of complex
structures modulo diffeomorphism.
(Such families do not in general exist, unless the moduli problem
represented by $\Mc$ has been formulated to
include level structures and polarizations, so this use of $\X_t$ is
strictly speaking an abuse of notation.)
There is a  natural diagram
$$\CD
\Ms @<<< \MN \\
 @.      @VVV \\
     @.  \Mc
\endCD$$
of ``forgetful'' maps which relates these moduli spaces.

The fibers of the map  $\MN\to\Ms$ are determined by the precise nature
of the holonomy group.
We henceforth restrict our attention to the case $h^{2,0}(\X_t)=0$, in which
the map $\MN\to\Ms$ is known to be finite (cf.\ \cite{\StringKThree}).
We still need to understand the
fibers of the map $\MN\to\Mc$.  To this end, let us fix a complex
structure $t$.  Then each $g_{ij}$
corresponding to a point in the fiber over $t$ determines a
{\em K\"ahler form}\/
$\omega:=(\sqrt{-1}/2)\sum g_{\alpha\bar\beta}
\,dz^\alpha\wedge d\bar z^\beta.$
This is a closed, nondegenerate, real $2$-form which can be regarded
as specifying a symplectic structure on $X$.
The de Rham classes of all possible $\omega$'s (called {\em K\"ahler
classes}) form an open convex
cone $\K_t\subset H^2_{\DR}(X,\R)$,
the {\em K\"ahler cone}\/  of $\X_t$.
  There is also a closely related cone $(\K_t)_+$, the
{\em nef cone}, defined by
$$(\K_t)_+:=\operatorname{Hull}\left(\overline{\K_t}\cap
H^2_{\DR}(X,\Q)\right).$$
(This cone includes the rationally defined subsets
of the boundary of $\overline{\K_t}$, while omitting any irrational
parts of the boundary.)

By the theorems of Calabi \cite{\Calabi} and of Yau \cite{\Yau},
each class $J$ in $\K_t$ uniquely determines
a metric $g_{ij}$ on $\X_t$ with holonomy in $\SU(n)$ whose associated
K\"ahler form $\omega$ lies in the class $J$.
The fibers of the map $\MN\to\Mc$ can thus be written in the form
$\Gamma_t\backslash\D_t$, where $\D_t=H^2_{\DR}(\X_t,\R)+i\,\K_t$, and
$\Gamma_t=H^2_{\DR}(X,\Z)\semidirect\Aut(\X_t)$.
For if we are given
$(g_{ij},B)$  associated to a particular complex structure $t$,
the corresponding {\em complexified K\"ahler class}\/
$B+i\,J$ naturally takes values in $H^2_{\DR}(\X_t,\Z)\backslash \D_t$.
But we must also
mod out by $\Aut(\X_t)$, to take care of diffeomorphisms which preserve
the complex structure $t$.

The K\"ahler cone has now appeared.  Soon, we will need
to go ``beyond'' it.

\head 2.~Three-point functions \endhead

The physical theory also determines two trilinear maps called
{\em three-point functions}.
The first of these is known as the {\em B-model three-point function},
since it can be calculated using
Witten's ``B-model'' \cite{\Witten}---a close
relative of the original $\sigma$-model.
This ``three-point function'' can be regarded as
 a trilinear map among certain bundles
on the complex moduli space $\Mc$, defined by using a universal
family $\pi:\X\to\Mc$.
The   arguments of the three-point function are (local)
sections of certain bundles:
we take $\alpha\in\Gamma(\F^{n-p+1}/\F^{n-p+2})$,
$\beta\in\Gamma(R^1\pi_*\T_{\X/\M})$,
$\gamma\in\Gamma(\F^{p}/\F^{p+1})$, where the $\F^p$ are the Hodge bundles
for the family $\pi$, and $\T_{\X/\M}$ is the relative holomorphic
tangent bundle
of the family.  The B-model three-point  function is then defined to be:
$$\bra \alpha,\beta,\gamma\ket_{\text{B-model}}:=
\int_{\X_t}\alpha\wedge\nabla\!_\beta\,\gamma
\in\Gamma(\O_{\M}),$$
where $\nabla$ is the Gauss--Manin connection
and $\nabla\!_\beta$ is the directional derivative determined from
$\beta$ via the Kodaira--Spencer isomorphism.

The second three-point function is called the {\em A-model three-point
function}.  (We will only define this in the case $h^{2,0}(\X_t)=0$.)
The definition involves a bit of a technical digression.  It will be
formulated in terms of certain invariants, called the {\em Gromov--Witten
invariants}, which measure the rational curves on $\X_t$.
Verifying mathematically
that these invariants exist and have the properties expected by
the physicists is an area of intense study, with much recent
progress \cite{\Katz, \appendix, \Ruan,
\RuanTian, \LiuMcDuff, \KontsevichManin, \Kontsevich}.
We will give a heuristic description of these invariants,
following Witten's original discussion \cite{\tsm, \topgrav}, and use
them as if they had all of the expected properties (including in
particular the ``multiple cover formula'' of \cite{\AM}, which
we build into our definitions).

 For each class $\eta\in H_2(X,\Z)$, consider
the moduli space of maps
$$\IM'_\eta:=\{\text{generically injective
holomorphic maps } f:\P^1\to \X_t
\text{ with }
[f(\P^1)]=\eta\}.$$
A na\"{\i}ve dimension estimate suggests that $\dim(\IM'_\eta)=\dim(X)$, and
in fact by a theorem of McDuff \cite{\McDuff},
this is true provided that one deforms the complex
structure $t$ to a nearby (non-integrable) almost-complex structure on $X$.
In order to formulate our heuristic description, we pass to
such a nearby deformation.

For each point $P\in\P^1$, there is an evaluation map $e_P:\IM'_\eta\to \X_t$
given by $e_P(f)=f(P)$.  Then the Gromov--Witten invariants should be
(heuristically) defined as
$$
G'_\eta(A,B,C):=e_0^*(A)\cup e_1^*(B)\cup e_\infty^*(C)|_{[\IM'_\eta]},
\tag\GW$$
for $A\in H^{p-1,p-1}(\X_t)$, $B\in H^{1,1}(\X_t)$, and
$C\in H^{n-p,n-p}(\X_t)$.  (Note:  if $p=1$ or $p=n$, then the intersection
in the moduli
space cannot be made transverse and
$G'_\eta(A,B,C)=0$.)
The difficulty with this attempted
definition is that $\IM'_\eta$ is not compact,
so quite a bit of care must be used in trying to evaluate \eqref{\GW}.
Much of the recent work on these invariants has been based on
Gromov's compactification \cite{\Gromov} of $\IM'_\eta$, but other
compactifications have also been used.  It is not yet clear that all
proposed definitions agree.

We will assume that Gromov--Witten invariants can be defined somehow,
and use them to define the A-model three-point
functions in the following way.
Assume for simplicity that $H^2(X,\Z)$ and $H_2(X,\Z)$ have no torsion.
Let $\cC\subset\K_t$ be an open cone of the
form $\R_{>0}e^1+\dots+\R_{>0}e^r$ for some
 basis $e^1,\dots,e^r$ of $L:=H^2(X,\Z)$, and let
$e_1,\dots,e_r$ be the dual basis of $H_2(X,\Z)$.  We express elements
of $L_{\R}+i\,\cC$ in the form $\sum a_je^j$, and the coefficients
$a_j$ must satisfy
$\Im(a_j)>0$. The class of this element
in $L\backslash (L_{\R}+i\,\cC)$ can then be described by the quantities
$q_j:=\exp(2\pi ia_j)$, which provide coordinates on
$L\backslash (L_{\R}+i\,\cC)$.
Note that those coordinates are subject to the constraint
$0<|q_j|<1$.  In fact, the space $L\backslash (L_{\R}+i\,\cC)$ is isomorphic to
$(\Delta^*)^r$, where $\Delta^*$ is the punctured disk.  It admits
a natural partial compactification $(\Delta^*)^r\subset\Delta^r$
with a distinguished boundary point $0$.

The three-point function associated to this region
$L\backslash (L_{\R}+i\,\cC)\subset\Gamma_t\backslash \D_t$
is given by:
$$\bra A,B,C\ket_{\text{A-model}}:=A\cdot B\cdot C
+\sum_{0\ne\eta\in H_2(X,\Z)}\frac{q^\eta}{1-q^\eta}\,G'_\eta(A,B,C),
\tag\threepoint$$
where $q^\eta$ denotes $\prod (q_j)^{(\eta^j)}$ when $\eta=\sum\eta^je_j$.
All nonzero terms in this series have $\eta^j>0$.  However,
no convergence properties of the series are known, so at present, the
three-point function
 must be considered to take values in the formal power series
ring $\C[[q_1,\dots,q_n]]$, the completion of the coordinate
ring of $\Delta^r$ at the distinguished
boundary point $0$.  If we could learn
something about the radius of convergence
of this function, we would gain some information
about the domain within $\Ms$ in which the quantum field theory converges.

Notice that in the case of dimension 3 (i.e., $n=3$), the three-point
function involves the Gromov--Witten invariants only when all three of
$A$, $B$, and $C$ come from $H^{1,1}(\X_t)$ (i.e., when $p=2$).
In this case, each Gromov--Witten invariant simply
 counts the number of rational curves in the corresponding homology
class (with appropriate signs),
multiplied by the degrees of the class with respect to the
given divisors.  In other words,
$$\bra A,B,C\ket_{\text{A-model}}=A\cdot B\cdot C
+\sum_{\eta\ne0}\frac{q^\eta}{1-q^\eta}\,(A\cdot\eta)(B\cdot\eta)(C\cdot\eta)
\,\#(\IM'_\eta/\PSL(2,\C)).$$
As written, this formula only makes sense for generic almost-complex
structures (for which $\IM'_\eta/\PSL(2,\C)$ is a finite set);
 even for them, care must be taken when calculating the ``number''
of points of the set, as some of them should be counted with a minus
sign.  (See \cite{\twoparamtwo, \S8} for a discussion of this issue.)

In spite of an apparent dependence on
the complex structure $t$, in fact the Gromov--Witten
invariants are independent of $t$.\footnote{This is simply because
the definition as we have formulated it
uses a generic nearby almost-complex structure.
It would be desirable to have definitions of these
invariants purely within algebraic geometry---for such definitions,
the independence from $t$ will be more difficult to verify.}
It follows that the A-model three-point
function depends on the ``complexified K\"ahler'' parameters $q_j$
(which are local coordinates on $\Gamma_t\backslash \D_t$)
 but not on the complex structure parameter $t$.

The cure for the apparent dependence on the choice of
cone $\cC$---called a {\em framing}\/ in \cite{\compact}---is more
difficult.  The most ideal circumstances are represented by varieties
which satisfy the

\proclaim{Cone Conjecture \cite{\compact}}
There is a rational polyhedral cone
$$\Pi\subset H^2(X,\R)$$
the union of
whose
translates $\gamma(\Pi)$ by automorphisms $\gamma\in\Aut(\X_t)$ covers
the nef cone $(\K_t)_+$.
\endproclaim

\noindent
(This conjecture was originally made to ensure that the space
$\Gamma_t\backslash \D_t$ could be partially compactified using a construction
of Looijenga \cite{\Looijenga} related to
 the Satake--Baily--Borel
compactification.)
The polyhedron $\Pi$ in the cone conjecture can in fact be chosen to
be the closure of a fundamental domain for the $\Aut(\X_t)$-action, and such a
$\Pi$ can be subdivided into cones $\cC_j$  which are generated
by bases of $L$; these cones $\cC_j$ can then be used in the definition of
three-point functions as above.

The cone conjecture clearly holds whenever the nef cone is itself rational
polyhedral
(a not uncommon occurrence), and it has been verified in at least
one nontrivial
example \cite{\GrassiMorrison}.

As an alternative to the cone conjecture,
it is possible to interpret the symbols $q^\eta$ as belonging
to the group ring $\C[H_2(X,\Z)]$, and to interpret the
three-point functions \eqref{\threepoint} as taking values in a
certain formal completion of that ring.\footnote{I am grateful to
A.~Givental for this remark.}  This provides another way to give
an intrinsic meaning to Eq.~\eqref{\threepoint}, independent of
choices of bases and cones.

\head 3.~Predictions from mirror symmetry \endhead

Mirror symmetry \cite{\Dixon, \LVW, \CLS, \GP}
predicts that Calabi--Yau manifolds often come in {\em mirror pairs}\/
$(X,Y)$, related by the existence
an isomorphism between nonlinear
$\sigma$-models
on $X$ and $Y$ which permutes the data in a certain specified
way.\footnote{The predicted relationship between the Euler numbers of these
manifolds is $e(Y)=(-1)^{\dim X}e(X)$, which leads to the Euler number sign
ambiguity (such as the $e(X)=\pm6$ issue mentioned above) present
in the early searches for specific Calabi--Yau threefolds.}  In particular,
the r\^oles of the moduli spaces $\Mc$ and $\Gamma_t\backslash \D_t$ should be
exchanged when passing from $X$ to $Y$,
and the three-point functions of A-model and B-model type
should be reversed.

This prediction has one very puzzling aspect:
whenever the cone conjecture holds,
 $\Gamma_t\backslash \D_t$ is a bounded domain, covered
by a finite number of ``punctured polydisks''
$L\backslash (L_{\R}+i\,\cC)\cong(\Delta^*)^r$.
On the other hand, for each choice of polarization, the open subset
$\Mc^{\text{pol}}\subset\Mc$ of polarized complex structures
has the structure of a quasi-projective variety (by a
theorem of Viehweg \cite{\Viehweg}).  These properties
would appear at first sight
to be incompatible with a mirror symmetry isomorphism.

There are two potential resolutions to this puzzle.  It may be that, after
shrinking $\D_t$ to the natural domain of definition of the
physical theory $\D$, there is a much larger symmetry group $\Gamma$
which acts on $\D$ (representing identifications between conformal
field theories which are not visible as identifications between nonlinear
$\sigma$-models), in such a way that $\Gamma\backslash \D$ is quasi-projective.
Alternatively, it may be that the physical theory can be analytically
continued beyond $\D_t$, and that when one attaches
$\Gamma_t\backslash \D_t$
to other parameter spaces for other types of conformal field theory,
one obtains a quasi-projective variety.\footnote{One can also imagine
combinations of these two scenarios, in which an enlargement of
 $\D_t$ admits an
action by a larger group.}

Analytically continuing outside $\D_t$ would take us ``beyond the
K\"ahler cone''.  This is what in fact happens in the physical
theory, as has been demonstrated recently in \cite{\phases, \catp}.
We explain this in the next section.

\head 4.~Flops and $\sigma$-models \endhead

\subhead 4.1.~The boundary of the K\"ahler cone \endsubhead

For the remainder of this paper, we will specialize to
 the case of algebraic threefolds, where the techniques of Mori
theory (cf.~\cite{\CKM}) are available for studying the K\"ahler cone.
A detailed analysis of the cone for Calabi--Yau threefolds
has been carried out by Kawamata
\cite{\Kawamata} and Wilson \cite{\Wilson}.
The boundary of the closure of $\K_t$
has certain {\em rational walls}\/ of codimension 1: each
is the intersection of the nef cone $(\K_t)_+$ with a hyperplane
$\Gamma^\perp$, where
 $\Gamma\subset \X_t$ is a curve which is
 {\em log-extremal}\/ in the sense of Mori theory.  Any linear system
whose numerical class lies at an interior point of  a
rational wall $\Gamma^\perp$ determines a
``contraction mapping'' which contracts to points all
effective irreducible curves $\Gamma'$ which
are numerically equivalent to $\lambda \Gamma$ for some
real number $\lambda$.
The contraction mappings associated to rational walls come in
several varieties:
\roster
\item {\em flopping contractions}\/ contract a finite number of curves
to points,

\item {\em divisorial contractions}\/ contract divisors to subvarieties of
lower dimension, and

\item {\em Mori fibrations}\/ have images of lower dimension.
\endroster
(We use similar names for the rational walls, calling them
  flopping walls, divisorial walls,
and Mori-fibration walls.)
We will primarily focus on flopping walls, but will briefly return to the
other two types at the end of the paper.

\subhead 4.2.~Flops \endsubhead

Our initial concern is to understand the behavior of the K\"ahler metric
$g_{ij}$
as we allow the complexified K\"ahler class $B+i\,J$ to approach a flopping
wall.
The simplest kind of flop is centered on a curve $\Gamma\subset \X_t$ with
normal bundle $\O(-1)\oplus\O(-1)$.  The behavior of the
metric itself is not known;\footnote{However, Candelas and de la Ossa
\cite{\CdO} have given a very interesting
local analysis of metrics in this situation.} we will settle
for a (local) analysis of the behavior of the associated symplectic structure.

To analyze the symplectic structure,
we first give a description of this flop in terms of variable
symplectic reductions of a fixed $\C^*$-action, similar in spirit to
\cite{\GS} and \cite{\Thaddeus}.
We begin with $\C^4$ with coordinates $(w,x,y,z)$,
and consider the action of $\C^*$ on $\C^4$ given by
$$(w,x,y,z)\mapsto(sw,sx,s^{-1}y,s^{-1}z)$$
for $s\in\C^*$.
Fix a symplectic form
$$\omega:=
{\sqrt{-1}\over2}\left(dw\wedge d\bar w+dx\wedge d\bar x+dy\wedge d\bar y
+dz\wedge d\bar z\right)$$
on $\C^4$.  There is then a moment map $\mu:\C^4\to\R$
for the action of $\C^*$, given
by
$$\mu(w,x,y,z):=\frac12\,(|w|^2+|x|^2-|y|^2-|z|^2).$$
The fibers of the moment map are invariant under the maximal compact subgroup
$S^1\subset\C^*$.
One can then form the {\em symplectic reductions}\/
$\mu^{-1}(r)/S^1$, for various values $r$ in the image of the moment map.
The symplectic form $\omega$ induces a symplectic form $\omega_r$
on the set of
smooth points of the reduced space $\mu^{-1}(r)/S^1$.

This $\C^*$-action can also be studied directly in algebraic geometry.
By a variant of the Kirwan--Ness theorem \cite{\Kirwan, \Ness},
for each value of $r$ there is an algebraic set $\Sigma_r\subset\C^4$
such that
the symplectic reduction $\mu^{-1}(r)/S^1$ is isomorphic to
the geometric invariant theory quotient
$(\C^4-\Sigma_r)\catquot\C^*$.
The set $\Sigma_r$ can be characterized as the union of all $\C^*$-orbits
whose closures are disjoint from $\mu^{-1}(r)$.
(The usual version of the Kirwan--Ness theorem
applies to group actions on projective varieties;
for the quasi-projective version used here, see \cite{\Kirwan, p.~115}.)

When $r=0$, $\Sigma_r$ is empty and the geometric invariant theory
quotient $\C^4\catquot\C^*$ is simply the spectrum of the ring of invariants.
Computing that ring is not hard: it is generated by the polynomials
$$
A:=wy,\quad
B:=wz,\quad
C:=xy,\quad
D:=xz
$$
subject to the relation
$$AD-BC=0;$$
it follows that
$$\C^4\catquot\C^*=\Spec\C[A,B,C,D]/(AD-BC).\tag{\sing}$$
If $r<0$, then $\Sigma_r$ is the set $\{y=z=0\}$.  It is not difficult
to check that the quotient coincides with
the blowup of
$\Spec\C[A,B,C,D]/(AD-BC)$
along the ideal $A=B=0$.
If $r>0$, then $\Sigma_r$ is the set $\{w=x=0\}$.  Similarly,
the quotient this time is
the blowup of
\eqref{\sing}
along the ideal $A=C=0$.

Moving from $r>0$ to $r<0$ is thus geometrically described as the
familiar process of blowing down the original curve
$\Gamma\subset\X_t$ to
a point, and then blowing back up in a different way to obtain
a new curve $\widehat \Gamma$ in a birationally equivalent
space $\widehat\X_t$.  (At $r=0$, the natural geometric model is
the singular one $\overline{\X}_t$, in which $\Gamma$
has been contracted to a point $P\in\overline{\X}_t$.)
Since the proper transform map gives a natural identification between
$H^2(\X_t)$ and $H^2(\widehat\X_t)$, we can regard the K\"ahler cones
$\K_t$ and $\widehat\K_t$ as lying in the same space.  When this is
done, the homology classes of $\Gamma$ and $\widehat\Gamma$ are related
by $[\widehat\Gamma]=-[\Gamma]$.  Moreover, by the
Duistermaat--Heckman theorem \cite{\DH}, the cohomology class of
$\omega_r$ moves along a piecewise linear path in the cohomology
group (cf.~\cite{\GS}),
with $[\omega_r]\in\K_t$ for $r>0$ and $[\omega_r]\in\widehat\K_t$
for $r<0$.\footnote{We are blurring the distinction between the form $\omega_r$
(constructed using the local analysis
near $\Gamma$) and global properties of its cohomology class $[\omega_r]$.
This does not affect our statements, however, thanks to the ``equivalence
between algebraic and $\sigma$-model coordinates'' discussed in \cite{\small}.}
  We thus see that the nef cones $(\K_t)_+$ and $(\widehat\K_t)_+$
meet along a common rational wall---itself naturally associated
to the singular space $\overline{\X}_t$---which is a flopping wall
for both.

\subhead 4.3.~The effect on three-point functions \endsubhead

Having found that the K\"ahler cones of $\X_t$ and $\widehat \X_t$ can be fit
together in a natural way, we  now examine the A-model three-point functions
in these two cases.
We assume that $\X_t$ is sufficiently generic so that all effective curves
in the numerical equivalence class $\Gamma$ have normal bundle
$\O(-1)\oplus\O(-1)$, and that each is disjoint from all other rational
curves on $\X_t$.  It then follows that
the birational correspondence between $\X_t$ and
$\widehat \X_t$ does not disturb any curves other than those numerically
equivalent to $\Gamma$ and $\widehat \Gamma$.  Let
$n_\Gamma=n_{\widehat \Gamma}$
be the number of such curves.  We can write the A-model
three-point functions of $\X_t$
in the form
$$\align
\bra A,B,C\ket_{\text{A-model}}^{\X_t}
=A\cdot B\cdot C
&+\frac{q^{[\Gamma]}}{1-q^{[\Gamma]}}\,(A\cdot\Gamma)(B\cdot\Gamma)
(C\cdot\Gamma)\,n_\Gamma \\
&+\sum\Sb \eta\in H_2(X,\Z)\\ \eta\ne\lambda\Gamma \endSb
\frac{q^\eta}{1-q^\eta}\,G'_\eta(A,B,C),
\tag\intheform
\endalign $$
and there is a similar expression for the three-point function of
$\widehat\X_t$, in which only the first two terms are different.
Note that $\D_t\subset\{q : |q^{[\Gamma]}|<1\}$ so that the first two
terms in the sum \eqref{\intheform} are well-defined throughout $\D_t$.
These first two terms  can clearly be analytically continued
to the region of $q$'s satisfying $|q^{[\Gamma]}|>1$, where they can
be compared with the first terms of the corresponding three-point
functions for $\widehat \X_t$.

\iffigs
\topinsert
$$
\matrix\epsfxsize=2in\epsfbox{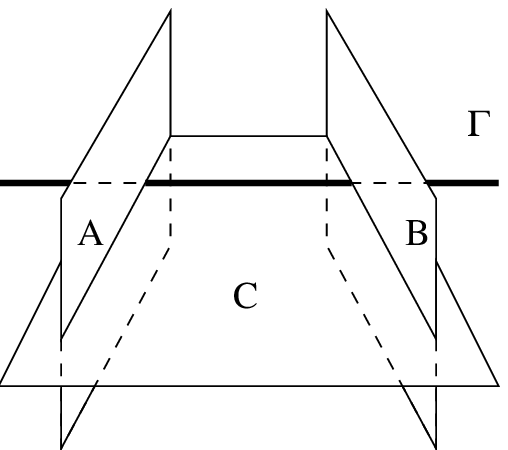} & \qquad &
\epsfxsize=2in\epsfbox{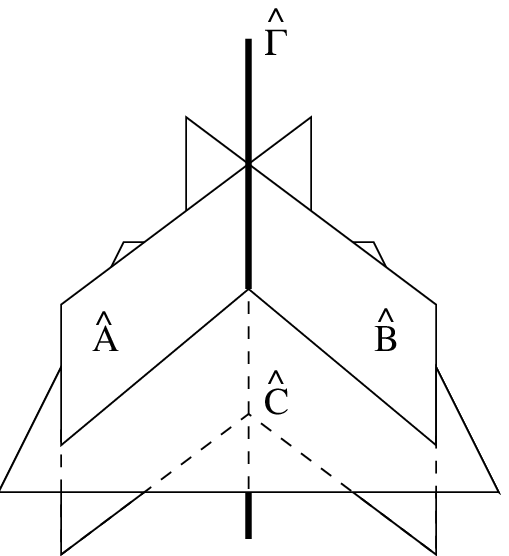} \cr
\quad & & \cr
\hbox{Figure 1a. Before the flop.} & &
\hbox{Figure 1b. After the flop.} \cr\endmatrix
$$
\endinsert
\fi

\proclaim{Lemma \cite{\phases, \small}}
For all $q\in(\D_t)_+\cup(\widehat\D_t)_+$ with $q^{[\Gamma]}\ne1$ we have
$$\align
A\cdot B\cdot C
&+\frac{q^{[\Gamma]}}{1-q^{[\Gamma]}}\,(A\cdot\Gamma)(B\cdot\Gamma)
(C\cdot\Gamma)\,n_\Gamma \\
= & {\widehat A}\cdot {\widehat B}\cdot {\widehat C}
+\frac{q^{[\widehat \Gamma]}}{1-q^{[\widehat \Gamma]}}\,
(\widehat A\cdot\widehat \Gamma)(\widehat B\cdot\widehat \Gamma)
(\widehat C\cdot\widehat \Gamma)\,n_{\widehat \Gamma},
\endalign$$
where $\widehat A$, $\widehat B$, and $\widehat C$ are the proper transforms
of $A$, $B$, and $C$.

In other words, the change in the topological term is precisely
compensated for by the change in the $q^{[\Gamma]}$ term.
\endproclaim

\demo{Proof}
First consider the case---illustrated in Fig.~1---in which $A$ and $B$ meet
$\Gamma$ transversally
(at, say, $a$ and $b$ points, respectively),
and $\widehat C$ meets $\widehat\Gamma$ transversally (at, say, $c$ points),
so that $C$
contains $\Gamma$ with multiplicity $c$.
(In Fig.~1 we illustrate the configuration of divisors in
  the case $a=b=c=1$.)
$A$ and $B$ have no intersection points along $\Gamma$, but both
$\widehat A$ and $\widehat B$ contain $\widehat \Gamma$, and they
meet $\widehat C$.  The total number of intersection points
of $\widehat A$, $\widehat B$ and $\widehat C$ (counted
with multiplicity) which lie in $\widehat \Gamma$ is thus $abc$.

Since a similar thing happens for each curve in the numerical equivalence
class, we see that
$$\widehat A\cdot\widehat B\cdot\widehat C
-A\cdot B\cdot C = abc\,n_\Gamma
=-(A\cdot\Gamma)(B\cdot\Gamma)(C\cdot\Gamma)\,n_\Gamma\tag\firstterms$$
(using $A\cdot\Gamma=a$, $B\cdot\Gamma=b$, $C\cdot\Gamma=-c$).
On the other hand, since $[\widehat\Gamma]=-[\Gamma]$ and
$n_{\widehat\Gamma}=n_\Gamma$, we can compute:
$$\align
\frac{q^{[\Gamma]}}{1-q^{[\Gamma]}}\,&(A\cdot\Gamma)(B\cdot\Gamma)
(C\cdot\Gamma)\,n_\Gamma
- \frac{q^{[\widehat \Gamma]}}{1-q^{[\widehat \Gamma]}}\,
(\widehat A\cdot\widehat \Gamma)(\widehat B\cdot\widehat \Gamma)
(\widehat C\cdot\widehat \Gamma)\,n_{\widehat \Gamma}\\
=&\frac{q^{[\Gamma]}}{1-q^{[\Gamma]}}\,
(A\cdot\Gamma)(B\cdot\Gamma)
(C\cdot\Gamma)\,n_\Gamma
+\frac{q^{-[\Gamma]}}{1-q^{-[\Gamma]}}
(A\cdot\Gamma)(B\cdot\Gamma)
(C\cdot\Gamma)\,n_\Gamma \\
=&\left(\frac{q^{[\Gamma]}}{1-q^{[\Gamma]}}
+\frac{1}{q^{[\Gamma]}-1}\right)
(A\cdot\Gamma)(B\cdot\Gamma)
(C\cdot\Gamma)\,n_\Gamma                    \tag\balance \\
=&-(A\cdot\Gamma)(B\cdot\Gamma)
(C\cdot\Gamma)\,n_\Gamma.
\endalign$$
Equating \eqref{\firstterms} and \eqref{\balance}
 proves the formula in this case.

To prove the formula in general, note it is linear in $A$, $B$, and $C$,
and so it suffices to prove the formula when
$A$, $B$, and $C$ are very ample divisors.  In particular we may assume
from the outset that $A$ and $B$ meet $\Gamma$ transversally, and that
$C\cdot\Gamma>0$.  Then $(-\widehat C)\cdot\widehat\Gamma>0$ as well,
and we can find ample divisors whose difference is $-\widehat C$.  Passing to
a multiple, we can in fact assume that there are very ample divisors
$\widehat H$ and $\widehat H'$ which meet $\widehat\Gamma$ transversally,
with $-\widehat C=\frac1N(\widehat H-\widehat H')$.  Applying the
formula for $(A,B,H)$ and for $(A,B,H')$ (which satisfy the hypotheses
of the special case), we deduce it for $(A,B,C)$.
\qed\enddemo

The conclusion to draw from this lemma is that {\em if}\/ the other
terms in the three-point function converge in some domain $\D$ contained
within $(\D_t)_+\cup(\widehat\D_t)_+$, then the entire three-point
function converges in $\D$, and it is the {\em same}\/ function for $\X_t$
as for $\widehat\X_t$.

\head 5.~Moduli spaces of birational $\sigma$-models  \endhead

We would like to construct a  moduli space
which incorporates this phenomenon of analytic continuation between
$\sigma$-models on birationally equivalent Calabi--Yau threefolds.
There are several technical difficulties in doing so, and the
discussion we give here can only be regarded as  a very preliminary attempt.
Ideally, the moduli problem would be formulated in terms of
birational metrics (or birational symplectic forms and birational
complex structures), modulo birational diffeomorphism.  Since we
don't quite understand how to do that, we will take a somewhat less natural
approach.

The first step is to fix a complex structure $t$ and construct
an enlargement of the space of complexified K\"ahler classes.
Our construction is based on a result of Kawamata \cite{\Kawamata},
who showed that the set of K\"ahler cones of all birational models
of $\X_t$, when transported to $H^2(\X_t,\R)$ via proper transforms
by all possible birational
maps, gives a chamber structure to the convex hull of their union.
(Taking the
convex hull simply adds the walls between adjacent chambers to the
union of the chambers.)
The resulting cone $\Mov(\X_t)$, called the {\em movable cone}\/
 by Kawamata,
can also be characterized as the interior of the
closely related cone
$$\Mov(\X_t)_+:=
\operatorname{Hull}\left(\overline{\Mov(\X_t)}\cap H^2_{\DR}(X,\Q)\right),$$
which is the
convex hull of the set of ``movable'' divisor classes---those
whose associated linear system has base locus of codimension at least
two.  The cone $\Mov(\X_t)_+$ can also be described as being
 the union of the proper transforms of the nef cones of all
birational models of $\X_t$.

We use the movable cone to define a {\em birational K\"ahler
moduli space}\/ of $\X_t$ in the form
$\Gamma_t^{\text{birat}}
\backslash\D_t^{\text{birat}}$,
where
$$\D_t^{\text{birat}}=H^2_{\DR}(\X_t,\R)+i\,\Mov(\X_t),$$
consists of all complex second cohomology
classes whose imaginary part lies in the movable
cone,
and where
$$\Gamma_t^{\text{birat}}=H^2_{\DR}(X,\Z)\semidirect\Bir(\X_t),$$
includes the entire group
 $\Bir(\X_t)$  of birational automorphisms of $\X_t$.
(Note that the action of this group will permute the various
chambers in Kawamata's chamber structure, that is, the various
K\"ahler cones of birational models of $\X_t$.)  As in the
case of the $\sigma$-model moduli space itself, we should regard
this space as only a first approximation to the true moduli space
of the physical theories.

In order for this space to be well-behaved, we need some control
over the action of the birational automorphism group.
As a natural generalization of both the cone conjecture and a
conjecture of Batyrev \cite{\BatyrevDuke}, we make the following

\proclaim{Birational Cone Conjecture}
There is a rational polyhedral cone
$${\cal P}\subset H^2(X,\R)$$
the union of whose
translates $\gamma({\cal P})$ by birational
automorphisms $\gamma\in\Bir(\X_t)$ covers
the cone $\Mov(\X_t)_+$.
\endproclaim

This is known to hold in at least one non-trivial
example---Calabi--Yau threefolds which are
fiber products of generic rational elliptic surfaces with section (as
studied by Schoen \cite{\schoen}).
The finiteness of the action of the birational automorphism group
on the set of birational models
was checked by Namikawa \cite{\Namikawa}, and the cone conjecture
was checked by Grassi
and the author \cite{\GrassiMorrison}.  Combining the two shows that
the birational cone conjecture holds in this case.

The lemma from section 4 shows that the A-model three-point functions
will generally have poles somewhere within
$\Gamma_t^{\text{birat}}
\backslash\D_t^{\text{birat}}$.
Points at which poles occur
 should be removed from
$\Gamma_t^{\text{birat}}
\backslash\D_t^{\text{birat}}$
if we hope to describe an actual moduli space
for the physical theories,
which would only be expected to contain points represented by
smooth Calabi--Yau manifolds.
On the other hand, to understand the natural limit points of this
moduli space we should leave
$\Gamma_t^{\text{birat}}
\backslash\D_t^{\text{birat}}$
intact
and try to construct a compactification (or at least a partial
compactification) of it.  This can be done using
 Looijenga's semi-toric
construction \cite{\Looijenga}
whenever the birational cone conjecture holds for $\X_t$.

The entire moduli space of birational $\sigma$-models should
somehow be constructed as
the union of the birational K\"ahler
moduli spaces:
$$\Mbir:=\bigcup_{t\in\Mc}
\Gamma_t^{\text{birat}}
\backslash\D_t^{\text{birat}}.$$
Unfortunately, at present, we do not know how either
$\Mov(\X_t)$ or $\Bir(\X_t)$
vary with parameters, so it is difficult
to topologize this space $\Mbir$ or discuss its properties in
any detail.

\head 6.~And beyond \endhead

We have seen that it is natural to go beyond the K\"ahler cone
of $\X_t$ to the movable cone, which describes physical theories
based on all birational models of $\X_t$.  Can we go even further?

In the case of K3 surfaces, the {\em elementary transformations}\/
play an important r\^ole in understanding the structure of the period
map \cite{\BurnsRapoport, \remarksKthree}.  These have a natural analogue
for Calabi--Yau threefolds \cite{\Wilson, \twoparam}, which lead
a bit beyond the movable cone.  Any divisorial contraction which
contracts a divisor $E$ to a curve $C$ of genus $g\ge1$ has
its associated divisorial wall of the form
$\Gamma^\perp$ with $\Gamma$ a generic
fiber of the induced map $E\to C$.  Since $E\cdot\Gamma=-2$,
there is an associated
{\em reflection}\/ in cohomology
$$H\mapsto H+(H\cdot\Gamma)E,$$
which can be used to reflect the movable cone of $\X_t$ through the wall
$\Gamma^\perp$.
As shown in \cite{\twoparam, \S9}, the A-model three-point functions
have a natural analytic continuation into this reflected cone,
compatible with the reflection mapping, and the physical theory
on the other side of the wall is isomorphic to the $\sigma$-model
theory on some birational model of $\X_t$.  We can thus extend
our moduli space to the ``reflected movable cone'', which includes
the images of $\Mov(\X_t)$ under all such reflections (on all
birational models of $\X_t$).

What about other walls?  For those Calabi--Yau threefolds which can be
realized as hypersurfaces in toric varieties, it is possible to
study the full
analytic continuation\footnote{More precisely, what is studied is the
full analytic continuation of that part of the K\"ahler moduli space
which comes from the ambient toric variety.} of the K\"ahler moduli space
\cite{\catp, \small}
by combining a special formulation of the physical theories on
those spaces due to Witten \cite{\phases}
with an analysis based on mirror symmetry.
When this is done, it is found that there are indeed analytic
continuations beyond other kinds of walls---in fact, the full moduli space
exhibits a vast chamber structure which includes many regions that meet the
original K\"ahler and movable cones only at the origin.  The currently
available
descriptions of the physical theories from
these other regions depend on the special formulation of \cite{\phases},
and appear to be linked to the ambient toric variety.  An abstract
description of these theories---or at least of their associated
cones---will be needed before any attempt
can be made to formulate a completely general account
of going beyond the ``reflected movable cone'' for arbitrary
Calabi--Yau threefolds.

\enddocument